\documentstyle[twoside,12pt,a4]{article}
\newfont{\Bb}{msbm10}
\newcommand{\CC}{\mbox{\Bb C}}

\newcommand{\ZZ}{\mbox{\Bb Z}}
\newcommand{\hol}{{\cal O}}
\newcommand{\HOL}{{\cal O}_X}

\newcommand{\DX}{{\cal D}_X}
\newcommand{\DER}{{{\cal D}{\rm er}_{\CC}(\HOL)}}

\newcommand{\D}{{\cal D}}

\newcommand{\vcero}{{{\cal V}_0^I}({\cal D})}
\newcommand{\VCERO} { {\cal V}_0^Y ({{\cal D}_X}) }

\newcommand{\DY}{{\cal D}_X[\star Y]}
\newcommand{\OY}{\HOL[\star Y]}

\newcommand{\derlog}{\mbox{\rm Der(log $f$)}}
\newcommand{\logf}{\log f}
\newcommand{\logY}{\log Y}

\newcommand{\derI}{\mbox{\rm Der(log $I$)}}
\newcommand{\derY}{\mbox{${\cal D}{\rm er}(\logY)$}}
\newcommand{\der}{{\rm Der}_{\CC}(\hol)}
\newcommand{\GrD}{\mbox{${\cal G}{\rm r}_{F^{\bullet}}({\cal D}_X)$}}
\newcommand{\grD}{\mbox{${\rm Gr}_{F^{\bullet}}({\cal D})$}}
\newcommand{\GRVo}{{\cal G}{\rm r}_{F^{\bullet}}\left({\cal V}_0^Y({\cal
D}_X)\right)}

\newcommand{\deriv}{\frac{\partial}{\partial x_1},
\frac{\partial}{\partial x_2},\cdots,\frac{\partial}{\partial x_n}}
\newcommand{\delsub}{{\delta_1},{\delta_2},\cdots,
   {\delta_n}}

\newcommand{\be}{{\begin{enumerate}}}

\newcommand{\ee}{{\end{enumerate}}}
\newcommand{\ec}{{\end{center}}}
\newcommand{\bc}{{\begin{center}}}

\newcommand{\g}{{\gamma}}
\newcommand{\de}{{\delta}}

\newcommand{\OM}{{\Omega}}

\newcommand{\V}{{{\cal V}}}

\newcommand{\bu}{\bullet}
\newcommand{\lo} {{logarithmic }}

\newcounter{numero}[subsection]
\renewcommand{\thenumero}{\thesubsection .\arabic{numero}}

\renewcommand{\thesubsection}{\arabic{section}.\arabic{subsection}}
\newcounter{numeron}[section]
\renewcommand{\thenumeron}{\thesection .\arabic{numeron}}
\renewcommand{\thesection}{\arabic{section}}
\newenvironment{corolario}{
\dimen255=\parindent \parindent=0in
\refstepcounter{numero}{\medskip}{\bf Corollary \thenumero.--}\
}{\parindent=\dimen255\par\vspace{1ex}}

\newenvironment{definicion}{
\dimen255=\parindent \parindent=0in
\refstepcounter{numero}{\medskip}{\bf Definition \thenumero.--}\
}{\parindent=\dimen255\par\vspace{1ex}}

\newenvironment{teorema}{
\dimen255=\parindent \parindent=0in
\refstepcounter{numero}{\medskip}{\bf Theorem \thenumero.--}\ }
{\parindent=\dimen255\par\vspace{1ex}}

\newenvironment{lema}{
\dimen255=\parindent \parindent=0in
\refstepcounter{numero}{\medskip}{\bf Lemma \thenumero.--}\ }
{\parindent=\dimen255\par\vspace{1ex}}

\newenvironment{proposicion}{
\dimen255=\parindent \parindent=0in
\refstepcounter{numero}{\medskip}{\bf Proposition \thenumero.--}\ }
{\parindent=\dimen255\par\vspace{1ex}}

\newenvironment{nota}{
\dimen255=\parindent \parindent=0in
\refstepcounter{numero}{\medskip}{\bf Remark \thenumero.--}\ }
{\parindent=\dimen255\par\vspace{1ex}}

\newenvironment{prueba}{{\bf Proof:}\ }{\hfill $\Box$\par\vspace{2ex}}

\title{Logarithmic differential operators and \\ logarithmic de Rham complexes 
  \\ relative to a free divisor \thanks{Supported by DGICYT PB94-1435} }

\author{Francisco J. Calder\'{o}n-Moreno}
\date{{\small Fac. Matem\'{a}ticas, Univ. de Sevilla, Ap 1160, 41080
Sevilla, Espa\~{n}a} \\
{\small E-mail: calderon@atlas.us.es}}

\begin{document}

\maketitle

{\Large {\bf Introduction}}
\vspace{0.5cm}
 
  \indent
    In the present work we prove a structure theorem for 
  operators of the 0-th term of the $\V^Y_{\bullet}$-filtration relative 
  to a free divisor $Y$ of a complex analytic variety $X$. As an application,
   we give a formula for the 
  logarithmic de Rham complex in terms of $\V_0^Y$-modules, 
  which generalizes the classical formula for the usual de Rham complex
 in terms of $\DX$-modules, and the formula of Esnault-Viehweg in the
  case that $Y$ is a normal crossing divisor. Using this, we give a sufficient 
  condition for perversity of the \lo de Rham complex. 
   Now we comment on the contents of each part of the paper: 
  
  In the first section, we recall the concepts of \lo derivation and
   \lo form, as well as free divisor, all of them due to Kyogi 
  Saito \cite{ksaito_log}, and the definition of the ring $\VCERO$ of
   \lo differential operators along $Y$.
  
  In the second part, we study the \lo operators in the case that $Y$ is free.
  We give a structure theorem in which we prove that the ring of \lo
  differential operators is the polynomial algebra generated by the 
  \lo derivations over the sheaf $\HOL$ of holomorphic functions. 
  As a consequence, $\VCERO$ is a 
  coherent sheaf. Thanks to this theorem, we can prove the equivalence between 
  $\VCERO$-modules and $\HOL$-modules with \lo connections. Therefore, an
  $\VCERO$-module (or \lo $\DX$-module) ${\cal M}$ defines a \lo de Rham 
  complex $\Omega^{\bullet}_X(\log Y)({\cal M})$.
  
 In the third part, we prove that the \lo 
 de Rham complex is canonically isomorphic to the 
 complex ${\bf R}{\cal H}{\rm om}_{\VCERO}(\HOL,{\cal M})$. To show this,
 we first construct a resolution of $\HOL$ as $\VCERO$-module, 
 which we call the \lo Spencer complex and denote by
 ${\cal S}p^{\bullet}{\rm (log}\ Y)$.
 
 Finally, we give a sufficient condition for perversity of the \lo de Rham 
 complex, 
which is a perverse sheaf if the symbols of a minimal generating
 set of \lo derivations form a regular sequence in
the graded ring associated
to the filtration by the order on ${\cal D}_X$. This condition always holds in
dimension 2.

 Some results of this paper have been announced 
 in \cite{nota-cras}. We give here the complete proofs of all of the 
 results announced in that note and other new results.
 
   {\it Acknowledgements:} I am grateful to David Mond for his interest 
 and encouragement. I wish to thank my advisor Luis Narv\'{a}ez for
 introducing
 me to the subject of this work and for giving
 me suggestions for the proofs of some of the results.

\section{Notations and Preliminaries}

  Let $X$ be a complex  analytic variety of dimension $n$, and $Y$ 
  a hypersurface of $X$ defined by the ideal ${\cal I}$. We will 
  denote by $\D_X$ the sheaf of linear differential operators over $X$,
  $\DER$ the sheaf of derivations of $\HOL$, and $\DY$ the sheaf 
  of meromorphic differential operators with poles along $Y$. Given a point 
  $x$ of $Y$, we will denote by $I=(f)$, $\hol$, $\der$ and $\D$ the respective 
  stalks at $x$. We will denote by $F^{\bu}$ the filtration of $\D_X$ 
  by the order of the operators and $\OM_X^{\bu}[\star Y]$ the
  meromorphic de Rham complex with poles along $Y$.

 \medskip 
 
\subsection{Logarithmic forms and logarithmic derivations. \\ 
 Free divisors}

 We are going to recall some notions of \cite{ksaito_log} that we will use 
 repeatedly:
 
 A section $\de$ of 
$\DER$, defined over an open set $U$ of $X$, is called a {\it \lo derivation }
(or vector field) if for each point $x$ in $Y\cap U$,
$\de_x({\cal I}_x)$ is contained in the ideal ${\cal I}_x$ (if 
$I={\cal I}_x=(f)$, it is sufficient that 
$\delta_x
(f)$ belongs to
$(f)\hol$).
 The sheaf of logarithmic derivations is denoted by $\derY$, and is a coherent
 $\HOL$-submodule of  $\DER$ and a Lie subalgebra.
We denote by $\derlog$, or $\derI$, the stalks at $x$ of $\derY$:
$$ \derlog =\{\delta\in {\rm Der}_{\CC}(\hol) \ / \ \delta(f) \in
(f)\}.$$

 We say that a meromorphic $q$-form $\omega$ with poles along $Y$, 
 defined in an open set  $U$,
is a \lo $q$-form  along  $Y$ or,
simply, a {\it \lo $q$-form}, if for every point $x$
in $U$, 
$f \omega$ and $df\wedge\omega$ are holomorphic at $x$.
The sheaf of \lo $q$-forms along $Y$ in $U$
is denoted by $\Omega^q_X({\rm log} \ Y)(U)$.
This definition gives rise to a coherent $\HOL$-module
$\Omega^q_X({\rm log}\ Y)$, whose stalks are:
$$\Omega^q({\rm log}\ f) = \Omega^q_X({\rm log}\
Y)_x=\{\omega\in \Omega^q_X[\star Y]_x \ / \ 
 f\omega\in \Omega^q,\
df\wedge\omega\in\Omega^{q+1}\}.$$

 The  \lo $q$-forms along $Y$ define
a subcomplex of the meromorphic de Rham complex along $Y$,
that we call the \lo de Rham complex and 
 denote by $\Omega_X^{\bullet}({\rm
log}\ Y)$.

 Contraction of forms by vector fields defines a perfect duality 
 between the $\HOL$-modules 
 $\Omega_X^1({\rm log}\ Y)$ and $\derY$, that we denote by 
 $\langle \ ,\ \rangle.$ Thus, both of them are reflexive.
 In particular, when $n={\rm dim}_{\CC}X=2$, $\Omega_X^1($log $Y$) and
\derY\ are locally free $\HOL$-modules of rank 2.

We say that $Y$ {\it is free at} $x$, or $I$ is a free ideal of
$\hol$, if \derI \ is free as $\hol$-module (of rank $n$).
 If
$f\in\hol$, we say that $f$ {\it is free} if the ideal $I=(f)$ is free.
We say that $Y$ is free if it is at every point $x$. In this case,
\derY \ is a locally free  $\HOL$-module of rank $n$.
 We can use the following criterion to determine when an hypersurface 
 $Y$ is free at $x$:

\medskip

{\bf Saito's Criterion:}
 The $\hol$-module \derlog \ is free if and only if there exist $n$
elements
$\delta_1,\delta_2,\cdots,\delta_n $ in \derlog,  with
$\delta_i=\sum_{j=1}^n a_{ij}(z)\frac{\partial}{\partial z_j} \
(i=1,\ldots,n),$ where $z = (z_1,z_2,\cdots,z_n)$ is a system of
coordinates of $X$ centered in $x$, such that the determinant
  $\det (a_{ij})$ is equal to $af$, with $a\in\hol$ a unit.
 Moreover, in this case,
$\{\delta_1,\delta_2,\cdots,\delta_n \}$ is a  basis 
of \derlog.
 
\smallskip

When $Y$ is free, we have the equality:
 $\Omega^p_X(\log Y) = \stackrel{p}{\wedge} \Omega^1_X(\logY)$. 
 Using the fact that $\Omega^1_X(\logY)\cong{\cal H}om_{\HOL}(\derY,\HOL),$
  we can construct a natural isomorphism: $$ \OM^p_X(\logY )\stackrel
{\g^p}{\cong}{\cal H}om_{\HOL}(\stackrel{p}{\wedge} \derY,\HOL), $$ 
\label{gamma} defined locally by 
$\gamma^p(\omega_1\wedge\cdots\wedge\omega_p)
(\delta_1\wedge\cdots\wedge\delta_p)=
\det (\langle\omega_i,\delta_j\rangle)_{1\leq i,j\leq p}.$

\subsection{${\cal V}$-filtration}

We define the ${\cal V}$-filtration relative to $Y$ on $\DX$ as
in the smooth case
 (\cite{mal_83}, \cite{kas_83_vfil}):
$${\cal V}_k^Y({\cal D}_X)=\{P\in {\cal D}_X \ /\ P({\cal
I}^j)\subset {\cal I}^{j-k}, \forall j\in\ZZ\}, \quad  k\in\ZZ,$$
where 
${\cal I}^p=\HOL$ when $p$ is negative. Similarly,
${\cal V}_k^I({\cal D})=\{P\in\D\ /\ P(I^j)\subset
I^{j-k},\forall j\in\ZZ \},$
with $k$ an integer, and $I^p={\cal O}$ when $p\geq 0$.
In the case of $I=(f)$, we note
${\cal V}_k^f({\cal D})={\cal V}_k^I({\cal D}).$

\begin{definicion}
A \lo differential operator (or, simplify, a \lo operator) is a
 differential operator of degree 0 with respect to the
${\cal V}$-filtration.
\end{definicion}

We see that:
$$\derY=\DER \cap \VCERO = {\cal G}r_{F^{\bu}}^1\left({\cal V}_0^Y
({\cal D}_X)\right),$$
$$F^1(\VCERO)=\HOL\oplus\derY,$$
where the last expression is consequence of 
$F^1(\DX)=\HOL\oplus\DER.$

\begin{nota}  \label{simgrad}
 The inclusion $\derY  \subset\GRVo $
gives rise to a canonical graded morphism of graded algebras: \\
$$\kappa:\ {\cal S}{\rm ym}_{\HOL}\left(\derY\right) \ \longrightarrow \
{\cal G}{\rm r}_{F^{\bu}}\left({\cal V}_0^Y({\cal D}_X)\right).$$
Similarly, we have a canonical graded morphism of
 graded $\hol$-algebras: \\ 
$\kappa_x:\ {\rm Sym}_{\cal O}\left(\derI\right) \ \longrightarrow
\ {\rm Gr}_{F^{\bu}}\left({\cal V}_0^I({\cal D})\right),$
which is the stalk of $\kappa$ at $x$.

\end{nota}

\section{ Logarithmic operators relative to a free divisor}

\subsection{The Structure Theorem}

We denote  by $\{\ ,\ \}$ the Poisson bracket defined in the graded
ring $\grD$
(cf. \cite{pha_79}, \cite{god_69}). Given two polynomials $F,G$ in
$\grD=\hol[\xi_1,\cdots,\xi_n]$:
$$\{F,G\}=\sum_{i=1}^n \frac{\partial F}{\partial \xi_i}
\frac{\partial G}{\partial x_i} -
\sum_{i=1}^n \frac{\partial F}{\partial x_i}
\frac{\partial G}{\partial \xi_i}. $$

\begin{proposicion}
Let $f$ be free. Consider a minimal system of generators 
$\{\delsub\}$  of  $\derlog $.
Let $R_0$  be a polynomial in $\grD$, homogeneous of order $d$, 
and such that there exist other polynomials $R_k$ in $\grD$, with
$k=1,\cdots,d$, homogeneous of order
$d-k$ such that:
\begin{equation} \{R_k,f\}=fR_{k+1},\ (0\leq k<d)
\label{tercera} \end{equation}
 (we will say that $R_0$ verifies the property (\ref{tercera})
 for $R_1,R_2,\cdots,R_d$). Then there exist polynomials $H_{j}^k$
in $\grD$, homogeneous of order $d-k-1$,
with $j=1,\cdots,n$ and $k=1,\cdots,d-1$,
such that:
\begin{enumerate}
\item[a)]  $ R_k=\sum_{j=1}^nH_j^k\sigma(\delta_j)$, \ where
$\sigma(\delta_j)$ denotes the principal symbol of $\delta_j$.
\item[b)]  $\{H_j^k,f\}=fH_j^{k+1} \ (1\leq j \leq n, \ 0\leq k<d-1)$.
 This is the same as saying: $H^k_j$ verifies the property (\ref{tercera}) for
$H^{k+1}_j,\cdots,H^{d-1}_j$.
\end{enumerate}
\end{proposicion}
\begin{prueba}
Let  $A=(\alpha_i^j)$ be the square  matrix  whose rows 
are the coefficients of the basis \{$\delsub$\} of \derlog \ with
respect to the basis $\deriv$ of $\DER$:
$$\delta_j=\sum_{i=1}^n\alpha_i^j\frac{\partial}{\partial x_i}
={\underline{\alpha}}^j \bullet \underline{\partial}^t,$$
with $j=1,\cdots,n$,  where we write $\underline{\partial}$
instead of $\left(\frac{\partial}{\partial x_1},
\cdots, \frac{\partial}{\partial x_n}\right)$.
We consider the ring ${\cal O}_{2n}=
\CC \{x_1,\cdots,x_2,\xi_1,\cdots,\xi_n\}.$
Thanks to the  Saito's Criterion, we know that the set 
$$\{\delta_1,\cdots,\delta_n,\frac{\partial}{\partial \xi_1},
\cdots,\frac{\partial}{\partial \xi_n}\}$$ 
is a basis of the
${\cal O}_{2n}$-module ${\rm Der}_{{\cal
O}_{2n}}$(log $f$).
So, as we have, for $k=1,\cdots,d$,
$$(f)\ni\{R_k,f\}=\sum_{i=1}^n
(R_k)_{\xi_i} f_{x_i},$$
 where
$f_{x_i}$ represents $\frac{\partial f}{\partial x_i}$ and
$(R_k)_{\xi_i}$ represents  $\frac{\partial R_k}{\partial \xi_i}$,
then there exist homogeneous
 polynomials  $G_j^k$ in $\grD$, of degree
$d-k-1$, or null, with $j=1,\cdots,n$ and $k=1,\cdots,
d-1$, such that
$$ \left((R_k)_{\xi_1},(R_k)_{\xi_2},\cdots,(R_k)_{\xi_n}\right)=
\sum_{j=1}^n G_j^k{\underline{\alpha}}^j.$$
Using the Euler relation
$R_k=\frac{1}{d}\sum_{i=1}^n(R_k)_{\xi_i}\xi_i$, and as
$\sigma(\delta_i)={\underline{\alpha}}^i\bullet\underline{\xi}^t $,
we obtain $$
 R_k=\frac{1}{d} \sum_{i=1}^n\sum_{j=1}^n G_j^k\alpha_i^j\xi_i
= \frac{1}{d}\sum_{j=1}^n G_j^k\sigma(\delta_j). $$
 By Saito's Criterion, the determinant of the  matrix $A$ is equal
to $uf$, with $u\in\hol$ invertible.
 Let $B=(b_{ij})=Adj(A)^t$.
 We have:
$$\left((R_k)_{\xi_1},(R_k)_{\xi_2},\cdots,(R_k)_{\xi_n}\right)=
\left(G_1^k,G_2^k,\cdots,G_n^k\right)A,$$ so
$$\left((R_k)_{\xi_1},(R_k)_{\xi_2},\cdots,(R_k)_{\xi_n}\right)B=
g\left(G_1^k,G_2^k,\cdots,G_n^k\right).$$

\noindent  Now:
$$ g\{G_j^k,f\}=\{gG_j^k,f\}= \sum_{i=1}^n f_{x_i}
\frac{\partial
(gG_j^k)}{\partial \xi_i}  = \sum_{i=1}^n  f_{x_i}
\sum_{l=1}^n
\frac{\partial
(R_k)_{\xi_l}}{\partial \xi_i} b_{lj} = $$ $$ \sum_{l=1}^n b_{lj}
\sum_{i=1}^n
\frac{\partial^2 R_k}{\partial \xi_l \partial \xi_i} f_{x_i} =
 \sum_{l=1}^n b_{lj}
\frac{\partial (\{R_k,f\})}{\partial \xi_l} =
 f \sum_{l=1}^n b_{lj}
\frac{\partial R_{k+1}}{\partial \xi_l} =
f\sum_{l=1}^n b_{lj} (R_{k+1})_{\xi_l}= $$
$$ f \sum_{l=1}^n b_{lj}
\sum_{p=1}^n G^{k+1}_p\alpha_l^p = 
 f \sum_{p=1}^n G_p^{k+1}
\sum_{l=1}^n b_{lj}\alpha_l^p = fgG_j^{k+1}.$$
\noindent
Therefore, $$\{G_j^k,f\}=fG_j^{k+1},$$ with $k=0,\cdots,d-2$ and 
$j=0,\cdots,n$.
 We conclude by setting
$H_j^k=\frac{1}{d}G^k_j$, for $j=1,\cdots,n$ and $k=0,\cdots,d-1$.
\end{prueba}
\begin{proposicion}
 Let be $\{\delsub \}$ 
 a basis of \derlog.  If a
 polynomial $R_0$ of $\grD$
 is homogeneous and verifies the property (\ref{tercera}) of the last 
 proposition, we can find a differential operator $Q$ in
 $\hol[\delsub]$ such that
 $R_0$ is the symbol of  $Q$.
\end{proposicion}

\begin{prueba}
  We will do the proof by induction on the order of $R_0$.
 If $R_0\in\hol$, it is obvious.
We suppose that the result holds if the order of $R_0$ is less than $d$.
 Now let $R_0$ of order $d$ verifying (\ref{tercera}).
By the last proposition there exist $n$ homogeneous polynomials $H^0_j$ 
of order
$d-1$ such that:
$$R_0=\sum_{j=1}^nH^0_j\sigma(\delta_j) , \ H_j^0\ {\rm verifies} \
({\rm \ref{tercera}}) \
(j=1,\ldots,n).$$
By induction hypothesis, there exist $Q_j \in \hol[\delsub]$ such that
 $ H_j^0= \sigma(Q_j)$. So
$$ R_0 =\sum_{i=1}^n\sigma(Q_i)\sigma(\delta_i) =
  \sum_{i=1}^n\sigma(Q_i\delta_i) =
\sigma(\sum_{i=1}^nQ_i \delta_i) = \sigma(Q) $$
and $ Q=\sum_{i=1}^n Q_i\delta_i \in \hol[\delsub]$.
\end{prueba}
\begin{nota}
Really, the previous argument proves that if $R_0$ verifies
(\ref{tercera}), then
$R_0$ is a polynomial in
 $\hol[\sigma(\delta_1),\cdots,\sigma(\delta_n)]$.
\end{nota}

\begin{teorema}  \label{teorema}
If $f$ is free and $\{\delta_1,\delta_2,\cdots,\delta_n\}$ is a basis 
of the $\hol$-module \derlog, each \lo operator 
$P$ can be written in a unique way as a polynomial 
$$P=\sum\beta_{i_1 \cdots i_n}
\delta_1^{i_1}\delta_2^{i_2}\cdots\delta_n^{i_n},\quad  \beta_{i_1
\cdots i_n} \in \hol .$$
 In other words, the ring of \lo operators is the
$\hol$-subalgebra of
$\D$ generated by \lo derivations:
$$\vcero=
\hol[\delta_1,\delta_2,\cdots,\delta_n]=\hol[\derlog].$$
\end{teorema}     \vspace{0.2cm}

\begin{prueba}
The inclusion $\hol[\delsub]\subseteq\vcero$ is clear. We will 
prove the other inclusion by induction on the order of $P_0 \in
\vcero$.  If the order of $P_0$ is zero, then it is a holomorphic 
function and the result is obvious.
 We suppose the result is true for every \lo operator 
$Q$ whose order is strictly less than $d$.
 Let $P_0$ be a \lo operator of order $d$.
 We know that:
$$[P_0,f]=fP_1,$$ with
$P_1\in\vcero$. So, there exist several $P_k$, with $k=0,\cdots,d$,
such that
$[P_k,f]=fP_{k+1}.$ If we set $R_k=\sigma(P_k)$, in the case that 
$P_k$ has order $d-k$,
and $R_k=0$ otherwise, we obtain:
$$\{R_k,f\}=\{\sigma_{d-k}(P_k),f\}=\sigma_{d-k-1}([P_k,f])=
f\sigma_{d-k-1} (P_{k+1})=fR_{k+1}.$$
 By the previous proposition, there exists
$Q$ in $\hol[\delsub]$ of order $d$ and such that $\sigma(P_0)=\sigma(Q)$.
As the order of $P_0-Q\in  \vcero$  is strictly less than $d$, we
apply the induction hypothesis to $P_0-Q$ and obtain
$$ P_0=P_0-Q+Q \in \hol[\delsub],$$
as we wanted. \\
On the other hand, using the structure of Lie algebra it is clear 
that we can write a \lo operator as a $\hol$-linear combination of the 
monomials $\{\delta_1^{i_1},\cdots,\delta_n^{i_n}\}$. The uniqueness of this
expression follows from the fact that 
these monomials are linearly independent over $\hol$.

\end{prueba}

\begin{nota}  \label{polinom}
 As a immediate consequence of the theorem (see the previous
 remark), we obtain an isomorphism:
$${\rm Gr}_{F^{\bu}}\left({\cal V}_0^I({\cal D})\right)
\stackrel{\alpha}{\cong}
\hol[\sigma(\delta_1),\cdots,\sigma(\delta_n)].$$
\end{nota}

\begin{corolario} \label{simgra}
If $Y$ is free at $x$, the morphism 
$ \kappa_x$ from the symmetric algebra ${\rm Sym}_{\cal O}(\derlog)$  to
${\rm Gr}_{F^{\bu}}\left({\cal V}_0^f({\cal D})\right)$
(see remark \ref{simgrad})
is an isomorphism of graded $\hol$-algebras.
As a consequence, if $Y$ is a free divisor,
the canonical morphism
$$\kappa: \quad {\cal S}{\rm ym}_{\HOL}\left(\derY\right) \to
{\cal G}r_{F^{\bu}}\left({\cal V}_0^Y({\cal D}_X)\right)$$
is an isomorphism.
\end{corolario}

\begin{prueba}
Let $x$ be in $X$ and $f\in\hol$ a local reduced equation of $Y$ at a
neighbourhood
of $x$. Let $\{\delta_1,\cdots,\delta_n \}$ be a basis of $\derlog$.
 $$\derlog=\oplus_{i=1}^n \hol \delta_i
\cong\oplus_{i=1}^n\hol\sigma(\delta_i).$$
 The symmetric algebra of the $\hol$-module
$\derlog$ is isomorphic to a polynomial ring: 
$${\rm Sym}_{\hol}\left(\derlog\right)\stackrel{\beta}{\cong}
\hol[\sigma(\delta_1),\cdots,\sigma(\delta_n)].$$
  We also have the
 inclusion:
$$\oplus_{i=1}^n\hol\sigma(\delta_i)
= {\rm Gr}^1_{F^{\bu}} \left(\vcero\right)\subset 
{\rm Gr}_{F^{\bu}}\left(\vcero\right),$$
where $\sigma(\delta_i)$ is the image of $\delta_i$ by the morphism
$\kappa_x$.
Therefore we conclude that the morphism $\kappa_x=\alpha^{-1}\beta$ is an 
isomorphism (see remark \ref{polinom}).
On the other hand, the inclusion
 $$ \derY = {\cal G}r_{F^{\bu}}^1\left({\cal
V}_0^Y
({\cal D}_X)\right) \subset\GRVo $$  gives rise to
a canonical graded morphism of graded $\HOL$-algebras 
(see remark \ref{simgrad}):
$\kappa:\ {\cal S}{\rm ym}_{\HOL}\left(\derY\right) \ \longrightarrow \
{\cal G}{\rm r}_{F^{\bu}}\left({\cal V}_0^Y({\cal D}_X)\right),$
whose stalk at each point $x$ of $Y$ is the canonical graded
isomorphism $\kappa_x$.
So, $\kappa$ is also an
isomorphism.
\end{prueba}

\begin{corolario} $\VCERO$ is a coherent sheaf of rings.
\end{corolario}
\begin{prueba}
By theorem 9.16 of \cite{bjo_79} (p. 83), we have only to 
prove that $\GRVo$
is coherent, but this sheaf is locally isomorphic to the 
polynomial ring  ${\cal O}_X[T_1,\cdots,T_n]$, 
which is coherent (\cite[lemma 3.2, VI, pg. 205]{ban_stan}).
\end{prueba}

\subsection{Equivalence between $\HOL$-modules with a \lo 
connection and left $\VCERO$-modules.}
\vspace{0.1cm}

\begin{definicion}  (cf. \cite{del_70})
Let ${\cal M}$ be a $\HOL$-module. A connection on
${\cal M}$, with \lo poles along $Y$, (or \lo connection on ${\cal M}$),
 is a $\CC$-homomorphism 
$\nabla$,
$$\nabla:\ {\cal M} \ \to \ \Omega^1_X({\rm log}\ Y)\otimes {\cal M},$$
that verifies Leibniz's identity:
$\nabla (hm)=dh\cdot m + h \cdot \nabla (m),$ where
$d$ is the exterior derivative over $\HOL$. We will note
$\Omega^q_X({\rm log}\ Y)({\cal M})  =
\Omega^q_X({\rm log}\ Y)\otimes {\cal M}$.
\end{definicion}
If $\delta$ is a \lo derivation along $Y$, it defines a
$\CC$-morphism:
$$
\begin{array}{ccc}
\derY &  \longrightarrow &  {\cal E}{\rm nd}_{\CC}({\cal M}),\\
\delta  & \mapsto           &
\nabla_{\delta} \\
\end{array}
$$
where $\nabla_{\delta}(m)= \langle\delta,\nabla(m)\rangle $

\begin{nota}
A \lo connection $\nabla$ on ${\cal M}$ gives rise 
to a morphism of $\HOL$-modules
$$\nabla^{\prime}:\ \derY \ \to \ {\cal H}{\rm
om}_{\CC}({\cal M},{\cal M}) $$ which verifies  
Leibniz's condition:
$\quad \nabla^{\prime}_{\delta}(fm)=\delta(f)\cdot m +  f\cdot
\nabla^{\prime}_{\delta}(m).$ \\
 Conversely, given 
$\nabla^{\prime}$ verifying this condition, we define
$$\nabla:{\cal M}\to\Omega^1_X({\rm log} \ Y)({\cal M}),$$
with $\nabla(m)$ the element of
$\Omega^1_X({\rm log} \ Y)({\cal M})={\cal H}{\rm
om}_{\HOL}\left(\derY,{\cal M}\right)$ such that:
$$\nabla(m)(\delta)=\nabla^{\prime}_{\delta}(m).$$
\end{nota}

\begin{definicion}
A \lo connection $\nabla$
is integrable
if, for each pair $\delta$ and $\delta^{\prime}$ of \lo derivations,
it verifies:
$$\nabla_{[\delta,\delta^{\prime}]}=[\nabla_{\delta},
\nabla_{\delta^{\prime}}],$$
where $[\ ,\ ]$ represents the Lie bracket in $\derY$
and the commutator in ${\cal H}{\rm
om}_{\CC}({\cal M},{\cal M})$.
\end{definicion}

Given a \lo connection $\nabla$ and the exterior derivative $d$,
we can construct a morphism:
$$\nabla^q:\Omega^q_X({\rm log}\ Y)({\cal M})\to
\Omega^{q+1}_X({\rm log}\ Y) ({\cal M}),$$ 
for each $q=1,\cdots,n$. If 
$\omega$ and $m$ are sections of the sheaves
 $\Omega^p_X({\rm log} \ Y)$ and
${\cal M}$:
$$\nabla^q(\omega\otimes m)=d\omega\otimes m + (-1)^q\omega\wedge
\nabla(m).$$
The integrability condition is equivalent to $\nabla^q\circ
\nabla^{q-1}=0$, for every $q$ (cf. \cite{del_70}).
\vspace{0.1cm}

\begin{definicion}
Let ${\cal M}$ be a ${\HOL}$-module, and $\nabla$ an integrable \lo
connection along $Y$ on ${\cal M}$.
 With the above notation, we call the \lo de Rham complex of 
 ${\cal M}$, and we denote by
$\Omega_X^{\bullet}({\rm log} \ Y)({\cal M})$, the complex (of
sheaves of $\CC$-vector spaces):
$$0\to {\cal M}\stackrel{\nabla}{\to}\Omega^1_X({\rm log} \ Y)
({\cal M})
\stackrel{\nabla^1}{\to} \cdots \stackrel{\nabla^{q-1}}{\to}
 \Omega^q_X({\rm log} \ Y)({\cal M})\stackrel{\nabla^q}{\to}  $$
$$ \Omega^{q+1}_X({\rm log} \ Y)({\cal M})
\stackrel{\nabla^{q+1}}{\to} \cdots
\stackrel{\nabla^{n-1}}{\to} \Omega_X^n({\rm log} \ Y)({\cal M})
\to 0.$$
\end{definicion}
\noindent
In the particular case where 
the $\HOL$-module
${\cal M}$ is equal to $\HOL$ and the \lo connection 
$\nabla$ is equal to the exterior derivative $d:\HOL
\to \Omega^1_X({\rm log} \ Y)$, the morphisms
$$\nabla^q:\Omega^q_X({\rm log} \ Y)\longrightarrow
\Omega^{q+1}_X({\rm log} \ Y),$$
define the \lo de Rham complex of Saito.\\

  We consider the rings 
  $R_0=\HOL\subset R_1$ and
  $R = \VCERO = \bigcup_{k\geq 0} R_k  \ (1\in R_0\subset R)$,
   with $R_k=F^k(\VCERO).$ The ring ${\cal G}r(R)$ is commutative and
verifies
   
  (1) The canonical morphism $\alpha:  Sym_{R_0}({\cal G}r^1(R)) \to 
   {\cal G}r(R)$, defined by   $\alpha(s_1\otimes\cdots\otimes s_t)=
  s_1\cdots s_t$, is an
  isomorphism (see Corollary 2.1.6). \\
   With these conditions, $R_1$ is an $(R_0,R_0)$-bimodule, and a Lie
algebra 
  ($[x,y] = xy-yx \in R_1$, because ${\cal G}r(R)$ is conmutative).  
  Moreover, $R_0$ is a sub-$(R_0,R_0)$-bimodule of $R_1$ such that the two
   induced structures of $R_0$-module over the
   quotient $R_1/R_0$ are the same. 
   
  Let ${\bf T}_{R_0}(R_1) =  R_0 \oplus R_1 \oplus (R_1 \otimes_{R_0} R_1)
   \oplus \cdots
    $ be  the tensor algebra of the $(R_0,R_0)$-bimodule $R_1$, and let
   $\psi: {\bf T}_{R_0}(R_1) \to  R$ be the canonical morphism defined by
the
    inclusion $R_1\subset R$. We  prove a reciprocal theorem of one 
    Poincar\'{e}-Birkhoff-Witt theorem 
    \cite[theorem 3.1,p.198]{rin_63} .
    
 \begin{proposicion}
   The morphism $\psi$ 
   induces an isomorphism:
  $$\phi: {\bf S}= \frac{{\bf T}_{R_0}\left( R_1 \right)}{J} 
\cong R, \quad  \phi ((i(x_1) \otimes \cdots \otimes
i(x_t))+J)=x_1x_2\cdots x_t,$$
 where $i$ the inclusion of $R_1$ in the tensor algebra, 
 and $ J $ is the two sided ideal generated by the elements:\\
  $$ {\rm a)}\ a-i(a),\ a\in R_0\subset R_1, \quad 
  {\rm b})\ i(x) \otimes i(y) - i(y) \otimes i(x) - i([x,y]), \  x,y\in
R_1.$$
 \end{proposicion}

 \begin{prueba}
  First, we check that the morphism 
  $\phi: {\bf S} \to R$ is well defined:
  \begin{center} $\psi(a-i(a)) = a-a = 0,\ a\in R_0,$\\
  $\psi ( i(x) \otimes i(y) - i(y) \otimes i(x) - i([x,y]) ) =
   xy-yx-[x,y] = 0, \ x,y\in R_1.$
    \end{center}
    The algebra ${\bf T}_{R_0}(R_1)$ is graded,
  so it is filtered, and induces a filtration on the quotient. The induced
  morphism 
  $\phi: {\bf S} \to R$
  is filtered:  $$\psi(a)=a\in R_0,\  
  \psi (i(x_1) \otimes \cdots \otimes i(x_t)) = x_1x_2\cdots x_t \in
R_t.$$ 
    So, we can define a graded morphism of $R_0$-rings.
 $$\pi: {\cal G}r \left( {\bf S} \right)  \to
   {\cal G}r(R),$$
   $$\pi(\sigma_t(i(x_1)\otimes \cdots \otimes i(x_t)+J)) =
    \sigma_t^{\prime}(x_1\cdots x_t) = 
    \overline{x_1}\cdots \overline{x_t},$$
   where $x_i\in R_1$, $\overline{x_i}=\sigma_1^{\prime}(x_1)$ is the
class of
    $x_i$ in $R_1/R_0$,
   $\sigma_t(P)$ is the
    class of 
   $P\in {\bf S}$ in ${\cal G}r^t({\bf S} )$,
    and $\sigma_t^{\prime}(Q)$ the class of 
   $Q\in R_t$ in ${\cal G}r^t(R)$. 
   Note that ${\cal G}r( {\bf S})$ is conmutative: 
  it is generated by the elements $\sigma_0(a+J)$, $\sigma_1(i(x)+J)$, 
  with $a\in R_0$,
  $x\in R_1$, and 
  \begin{center}   $[i(x)+J,i(y)+J]=i([x,y])+J$,\\
  $[a+J,i(x)+J]=i(ax-xa)+J=b+J,\ b=ax-xa\in R_0 .$ 
  \end{center}
  
    On the other hand, the image of $R_0\subset R_1$ in ${\bf S}$ 
  is exactly 
  the part 
  of degree zero of ${\bf S}$, and then we obtain a morphism of 
  $R_0$-modules
  from ${\cal G}r^1(R) = R_1/R_0$ to ${\cal G}r^1( {\bf S})$ which 
  induces a morphism of $R_0$-algebras:
  $$\rho:  {\cal S}ym_{R_0}\left(\frac{R_1}{R_0}\right) \to
   {\cal G}r\left( {\bf S} \right), $$
   $$ \rho(\overline{x_1}\otimes\cdots \otimes\overline{x_t}) = 
   \sigma_t(i(x_1)\otimes\cdots\otimes i(x_t)+ J),$$
   which is obviously surjective.
  The composition $\pi\rho$ is equal to $\alpha$, and, by property (1) of
$R$,  
 we deduce that $\rho$ is injective. 
    As $\rho$ and $\pi\rho$ are isomorphisms, $\pi$ is as well, 
  as we wanted to prove. 
\end{prueba}  
   
  \begin{corolario}
 Let  $Y$ be a free divisor. Let ${\cal M}$ be a $\HOL$-module.
  An integrable logarithmic connection on ${\cal M}$ gives rise to a left
   $\VCERO$-structure on ${\cal M}$, and vice versa.
 \end{corolario}

 \begin{prueba}
 A $\HOL$-module ${\cal M}$ with an integrable logarithmic connection
 $\nabla$
 has a natural structure of left $\VCERO$-module defined by
 its  structure as $\HOL$-module.
   Let $\mu$ be the morphism of $(\HOL,\HOL)$-bimodules:
   $$\mu: R_1=\HOL\oplus\derY \to {\cal E}{\rm nd}_{\CC}({\cal M}),\ \ 
   \mu(a)(m)=am, \ \ \mu(\delta)(m)=\nabla_{\delta}(m) .$$
  $\mu$ induces a morphism
  $\nu: {\bf T}_{R_0}(R_1) \to {\cal E}{\rm nd}_{\CC}({\cal M}),$
  and, as $\nu(J)=0$, we have a morphism
  $$\VCERO\simeq \frac{{\bf T}_{R_0}(R_1)}{J} \to 
  {\cal E}{\rm nd}_{\CC}({\cal M}),$$
  which defines an structure of $\VCERO$-module on ${\cal M}$.

   On the other hand, a left $\VCERO$-module structure on 
${\cal M}$ defines an 
integrable logarithmic connection 
$\nabla$ on the $\HOL$-module ${\cal M}$:
$$\nabla:\derY  \to  {\cal E}{\rm nd}_{\CC}({\cal M}), \quad
 \nabla_{\delta}(m)\ = \delta \cdot m.$$ 
\end{prueba}

\begin{nota} A left $\VCERO$-module structure on ${\cal M}$ defines 
a \lo de Rham complex.
In local coordinates
$(U;x_1,\cdots,x_n)$, with
$\{\delta_1,\cdots,\delta_n\}$ a local basis of $\derY$ and
$\{\omega_1,\cdots,\omega_n\}$ its dual basis, 
the differential of the complex is defined by:
$$\nabla^p(U)(\omega\otimes m)=d\omega\otimes m +
\sum_{i=1}^n \left((\omega_i\wedge\omega)\otimes
\delta_i\cdot m\right),$$
for any sections $\omega\in\Omega_X^1(\log Y) $ and
$m\in{\cal M}$.
In the particular case of the left
$\VCERO$-module $\HOL$, defined as $\VCERO$-module 
in a natural way
($P\cdot g=P(g),$
with $g$ a holomorphic function and $P$ a \lo operator), 
this canonical structure of $\HOL$ as left $\VCERO$-module 
is obviously equivalent to
 the integrable \lo connection over $\HOL$ defined naturally 
by the exterior derivative ($\nabla=d$):
$$\nabla_{\delta}(g)=\langle\delta,dg\rangle =\delta(g).
$$
\end{nota}

\section{The Logarithmic de Rham Complex}

 In this section, $Y$ will be a free divisor. 

\subsection{ The Logarithmic Spencer Complex} \label{Splog}

\begin{definicion} We call the \lo  Spencer complex, and denote by
${\cal S}p^{\bullet}($log $Y$), the complex:
$$0\to\VCERO\otimes_{{\cal
O}_X}\stackrel{n}{\wedge}\derY\stackrel{{\textstyle \varepsilon}_{-n}}
{\to}\cdots\quad \quad \quad \quad \quad \quad \quad \quad \quad \quad \quad$$
$$ \quad \quad \quad \quad \quad \quad \quad \quad \quad \quad \quad
\quad \quad \cdots \stackrel{{\textstyle \varepsilon}_{-2}}
{\to}\VCERO\otimes_{{\cal
O}_X}\stackrel{1}{\wedge}\derY\stackrel{{\textstyle \varepsilon}_{-1}}
{\to}
\VCERO,$$
 where $$
\varepsilon_{{\scriptscriptstyle
-p}}(P\otimes(\delta_1\wedge\cdots\wedge\delta_p))
=\sum_{i=1}^p
(-1)^{i-1} P\delta_i\otimes(\delta_1\wedge\cdots\wedge\widehat{\delta_i}
\wedge\cdots\wedge\delta_p) + $$
$$\sum_{1\leq i<j\leq p}
(-1)^{i+j}P\otimes([\delta_i,\delta_j]\wedge\delta_1\wedge\cdots\wedge
\widehat{\delta_i}\wedge\cdots\wedge\widehat{\delta_j}
       \wedge\cdots\wedge\delta_p), \ \ (2\leq p\leq n). $$
$$\varepsilon_{{\scriptscriptstyle -1}} (P\otimes\delta)=P\delta.$$
We can augment this complex of left $\VCERO$-modules by another morphism: 
$$\varepsilon_0:\VCERO\to \HOL,\ \  \varepsilon_0(P)=P(1).$$
We call the new complex $\widetilde{{\cal
S}}p^{\bullet}$(log $Y$).
\end{definicion}
 This definition is essentially the same as the definition of the
 usual Spencer complex ${\cal
S}p^{\bullet}$ of $\HOL$ (cf. \cite[2.1]{meb_formalisme})
and generalizes the definition given by
Esnault and  Viehweg {\cite[App.~A]{es_vi_86}} in the case of a
normal crossing divisor. We denote by
${\cal S}p^{\bullet}[\star Y]={\cal D}_X[\star
Y]\otimes_{{\cal D}_X} {\cal S}p^{\bullet}$ the meromorphic
Spencer complex of $\OY$.

\begin{teorema}          \label{splogres}
The complex ${\cal S}p^{\bullet}(\log Y$) is a
locally free resolution of
$\HOL$ as left $\VCERO$-module.
\end{teorema}
\begin{prueba}
To see the exactness of $\widetilde{{\cal S}}p^{\bullet}(\log Y$) 
we define 
a discrete filtration  $G^{\bullet}$ such that it induces an exact graded 
complex (cf. \cite[lemma 3.16]{bjo_79}):
$$G^k\left(\VCERO\otimes\stackrel{p}{\wedge}\derY\right)
=F^{k-p}\left(\VCERO\right)\otimes\stackrel{p}{\wedge}\derY,$$
$$G^k(\HOL)=\HOL.$$
We have
$${\cal
G}{\rm r}_{G^{\bu}}\left(\VCERO\otimes\stackrel{p}{\wedge}\derY\right)=
\GRVo[-p]\otimes\stackrel{p}{\wedge}\derY,$$
$${\cal G}{\rm r}_{G^{\bu}}(\HOL)=\HOL.$$
As the above filtrations are compatible with the differential 
of the complex $\widetilde{{\cal S}}p^{\bullet}(\log Y)$, 
we can consider the complex
${\cal G}{\rm r}_{G^{\bu}}\left(\widetilde{{\cal S}}p^{\bullet}(\log
Y)\right):$
$$0\to\GRVo[-n]\otimes_{{\cal
O}_X}\stackrel{n}{\wedge}\derY\stackrel{\psi_{-n}}{\to}\cdots$$
$$\stackrel{\psi_{-2}}{\to}\GRVo[-1]\otimes_{{\cal
O}_X}\stackrel{1}{\wedge}\derY\stackrel{\psi_{-1}}{\to}
\GRVo\stackrel{\psi_0}{\to}\HOL\to 0,$$
where the local expression of the differential is defined by:
$$\psi_{-p}(G\otimes\delta_{j_1}\wedge\cdots\wedge\delta_{j_p})=
\sum_{i=1}^p
(-1)^{i-1}
G\sigma(\delta_{j_i})\otimes\delta_{j_1}\wedge\cdots\wedge
\widehat{\delta_{j_i}}
\wedge\cdots\wedge\delta_{j_p},\ \  (2\leq p\leq n). $$
$$\psi_{-1}(G\otimes\delta_i)=G\sigma(\delta_i),\ \ \psi_0(G)=G_0,$$
with $ \{ \delta_1,\cdots,\delta_n\}$ a (local) basis of \derY.
This complex is the Koszul complex of the ring
$$\GRVo \cong {\cal S}{\rm ym}_{\HOL}\left(\derY\right) $$
 with respect  to the $\GRVo$-regular sequence
$\sigma(\delta_1),\cdots,\sigma(\delta_n)$ in the ring
$\GRVo$. Consequently, it is exact.
\end{prueba}

\begin{lema} \label{cambiarf}
 For every \lo operator $P\in {\cal
V}_0^f({\cal D})$, there exist, 
for each integer $p$, a \lo operator $Q\in
{\cal V}_0^f({\cal D})$ and an integer $k$ such that
$f^{-p} P = Q f^{-k}$.
\end{lema}
\begin{prueba}
We will prove the lemma by induction on the order of the \lo
operator. If $P$ has order 0, it is in $\hol$, and it is clear that
$f^{-p}P=Pf^{-p}$.  
Let $P$ be of order $d$, and 
consider the \lo operator $[P,f^p]$, of order $d-1$. By
induction hypothesis, there exists an integer $m$ such that:
$$[P,f^{-p}]f^m \in {\cal V}_0^f({\cal D}).$$
Let $k$ be the greatest of the integers $m$ and $p$. It is clear
that:
$$f^{-p}Pf^k = Pf^{k-p}-[P,f^{-p}]f^k \in {\cal V}_0^f({\cal D}).$$
This proves the result:
$Q=Pf^{k-p}-[P,f^{-p}]f^k$.
\end{prueba}

\begin{nota}
For every operator $Q$ in $\DY_x$, we can always find a strictly positive 
integer $m$ such that
 $f^mQ\in {\cal V}_0^f({\cal D})$. Equivalently, for each 
 meromorphic differential operator
 $Q$, there exists a positive integer $p$ and a \lo operator
 $Q^{\prime}$ such that we can write:
$$ Q=f^{-p} Q^{\prime}.$$
\end{nota}

Now we introduce several morphisms that we will use later.

\begin{lema}   \label{DXpY}  We have the following isomorphisms:

\begin{enumerate}

\item $\OY\otimes_{\HOL}\VCERO \stackrel{\sim}{\hookrightarrow} \DY
  \stackrel{\sim}{\hookleftarrow}
              \VCERO\otimes_{\HOL}\OY.$

\item  $\alpha:\DY\otimes_{\VCERO}
\HOL \cong \OY, \quad \quad \alpha(P\otimes g) = P(g).$

\item $\rho:\DY\otimes_{\VCERO}\DY\cong\DY,
\quad \quad \rho(P\otimes Q) = P Q.$

\end{enumerate}

\end{lema}

\begin{prueba}

1. The inclusions $\VCERO,\OY \subset \DY$ give rise to the previous
isomorphisms of $(\VCERO,\OY)$-modules. Locally:
$$a f^{-k}\otimes P = a f^{-k} P = a Q \otimes f^{-p} ,$$
with $P$ and $Q$ \lo operators such that $f^{-k}P = Q f^{-p}$.
We have seen how to obtain $Q$ from $P$ (lemma \ref{cambiarf}), and
we can obtain $P$ from $Q$ in the same way.
On the other hand, we saw in the previous remark how to express
a meromorphic differential operator as a product of a meromorphic
function and a \lo operator.

2. We have to compose the following
isomorphisms of left $\DY$-modules:
%$$ \DY\otimes_{\VCERO}\HOL \cong
$$\OY\otimes_{\HOL}\VCERO\otimes_{\VCERO}\HOL  \cong
\OY\otimes_{\HOL}\HOL \cong \OY. $$

3. We obtain this isomorphism of $\DY$-bimodules from the
 composition of the following isomorphisms:
%$$  \DY\otimes_{\VCERO}\DY \cong
 $$ \OY\otimes_{\HOL}\VCERO\otimes_{\VCERO}\DY \cong
\OY\otimes_{\HOL}\DY \cong $$
$$  \OY\otimes_{\HOL}\OY\otimes_{\HOL}\VCERO \cong
  \OY\otimes_{\HOL}\VCERO
\cong \DY, $$
where the isomorphism $\OY\otimes_{\HOL}\OY \cong \OY$  sends
(locally) the tensor product  $g_1\otimes g_2$ to the
meromorphic function  $g_1 g_2$.
\end{prueba}

\vspace{0.1cm}
\begin{proposicion} \label{Sp*Y}  We have the following isomorphisms of 
complexes of $\DY$-modules:  
\begin{enumerate}
 \item[(a)] $ {\cal D}_X[\star
Y]\otimes_{{\cal V}_0^Y({\cal D}_X)}
 {\cal S}p^{\bullet} \cong {\cal S}p^{\bullet}[\star Y].$
 \item[(b)] $ {\cal D}_X[\star Y]
\otimes_{{\cal V}_0^Y({\cal D}_X)}
{\cal S}p^{\bullet}(\logY)\cong  {\cal
S}p^{\bullet}[\star Y].$
 \end{enumerate}

 \end{proposicion}

\begin{prueba}
(a) As ${\cal S}p^{\bullet}$ is a subcomplex of ${\cal
D}_X$-modules of
${\cal S}p^{\bullet}[\star Y]$, and $\DY$ is flat over 
$\VCERO$, the complex $\DY\otimes_{\VCERO}{\cal
S}p^{\bullet}$ is a subcomplex of 
$\DY\otimes_{\VCERO}{\cal S}p^{\bullet}[\star Y],$ 
(see  lemma  \ref{DXpY}, 1.).
 But, by the third isomorphism of lemma \ref{DXpY}, this complex
is the same as ${\cal S}p^{\bullet}[\star Y]$. Hence, we have an 
injective morphism of complexes:
$$\DY\otimes_{\VCERO}{\cal S}p^{\bullet}  \longrightarrow
{\cal S}p^{\bullet}[\star Y],$$
defined locally in each degree by:
$P\otimes Q\otimes\delta_1\wedge\cdots\wedge\delta_ p \mapsto
PQ\otimes (\delta_1\wedge\cdots\wedge\delta_p).$
This morphism is clearly surjective and, consequently, an isomorphism.

(b) We consider $\VCERO$ as a subsheaf of $\hol$-modules of
${\cal D}_X$. Using the fact that
$\stackrel{p}{\wedge}\derY$ is $\HOL$-free, we have an
inclusion
$$\VCERO\otimes_{\HOL} \stackrel{p}{\wedge}\derY \ \hookrightarrow \
{\cal D}_X\otimes_{\HOL}\stackrel{p}{\wedge}\derY.$$
On the other hand, as $Y$ is free, we have a natural injective morphism from
$\stackrel{p}{\wedge}\derY$ to 
$\stackrel{p}{\wedge}\DER$
(cf. \cite[AIII 88, Cor.]{bou_ac_3_4}).
As ${\cal D}_X$ is flat over $\HOL$, we have other inclusion:
$${\cal D}_X\otimes_{\HOL}\stackrel{p}{\wedge}\derY \
\hookrightarrow
\ {\cal D}_X\otimes_{\HOL}\stackrel{p}{\wedge}\DER \ (p\geq 0).$$
Composing both of them, we obtain a new inclusion:
$$\VCERO\otimes_{\HOL}\stackrel{p}{\wedge}\derY \ \hookrightarrow
\
 {\cal D}_X\otimes_{\HOL}\stackrel{p}{\wedge}\DER,$$
for $p=0,\cdots,n$. These inclusions give rise to
an injective morphism of complexes of $\VCERO$-modules
$${\cal S}p^{\bullet}(\logY)\ \hookrightarrow \
{\cal S}p^{\bullet}.$$
As $\DY$ is flat over $\VCERO$ (see lemma \ref{DXpY}, 1.)
we have an injective morphism of complexes of $\DY$-modules:
$$\theta^{\prime}:\DY\otimes_{\VCERO}{\cal S}p^{\bullet}
(\logY)\ \hookrightarrow \ \DY\otimes_{\VCERO}{\cal
S}p^{\bullet},$$
defined by:
$\theta^{\prime}\left(P\otimes Q\otimes
(\delta_1\wedge\cdots\wedge\delta_p)\right)
= P\otimes Q\otimes (\delta_1\wedge\cdots\wedge\delta_p).$
This morphism is surjective, given $P$ local section of $\DY$, $Q$ in 
${\cal D}$ and $\delta_1,\cdots,\delta_n$ in $\der$, we have:
$$ P\otimes Q\otimes(\delta_1\wedge\cdots\wedge\delta_p)=\theta^{\prime}
\left( (Pf^{-k})\otimes Q^{\prime}\otimes (f\delta_1\wedge\cdots\wedge
f\delta_p)\right),$$
with $k>0$ and
$Q^{\prime}$ a local section of $\VCERO$ verifying
$f^k Q=Q^{\prime}f^p$ (see lemma \ref{cambiarf}).
Composing $\theta^{\prime}$ with the isomorphism of (a),
we obtain the isomorphism:
$$ \theta: \DY\otimes_{\VCERO}{\cal S}p^{\bullet}
(\logY)\ \stackrel{\sim}{\rightarrow} \
{\cal S}p^{\bullet}[\star Y],$$
with local expression:\
$ \theta (P\otimes Q\otimes (\delta_1\wedge\cdots\wedge\delta_p))
= P Q\otimes (\delta_1\wedge\cdots\wedge\delta_p).$
\end{prueba}
\vspace{0.1cm}

\subsection{The Logarithmic de Rham Complex}

For each divisor $Y$, we have a standard canonical isomorphism:
$${\cal H}{\rm om}_{\HOL}\left(\stackrel{p}{\wedge}\derY,{\cal
O}_X\right)
\stackrel{\lambda^p}{\cong}
{\cal H}{\rm om}_{\VCERO}
\left(\VCERO\otimes_{\HOL}\stackrel{p}{\wedge}\derY,\HOL\right),$$
defined by:
$\lambda^p(\alpha)(P\otimes\delta_1\wedge\cdots\wedge\delta_p) =
P\left(\alpha(\delta_1\wedge\cdots\wedge\delta_p)\right).$

Composing this isomorphism with the isomorphism $\gamma^p$ defined
in section \ref{gamma},
we can construct a natural morphism
$\psi^p=\lambda^p\circ\gamma^p:$
$$\Omega^p_X({\rm log}\ Y)\stackrel{\psi^p}{\cong}
{\cal H}{\rm om}_{\VCERO}
\left(\VCERO\otimes\stackrel{p}{\wedge}\derY,\HOL\right), \
$$
for $p=0,\cdots,n$. Locally:
$$\psi^p (\omega_1\wedge\cdots\wedge\omega_p)
(P\otimes\delta_1\wedge\cdots\wedge\delta_p) = P\left(
\det (\langle\omega_i,\delta_j\rangle)_{1\leq i,j
\leq p} \right). $$
with  $\omega_i\  (i=1,\cdots,n)$  local sections of
$\Omega^1_X({\rm log}\ Y)$ and $P$ a \lo  operator.

Similarly, if ${\cal M}$ is a left $\VCERO$-module, 
given an integer $p\in \{1,\cdots,n\}$, there exist the following
canonical isomorphisms:
$$ \gamma^p_{{\cal M}}: \ \Omega^p_X({\rm log}\ Y)\otimes_{\HOL}{\cal M}
\stackrel{\sim}{\to}
{\cal H}{\rm om}_{\HOL}\left(\stackrel{p}{\wedge}\derY,{\cal
M}_X\right),$$
%que viene dado localmente por
%$$\gamma^p_{{\cal M}}(\omega_1\wedge\cdots\wedge\omega_p\otimes m)
%(\delta_1\wedge
%\cdots\wedge\delta_p)=\det
%\left(\langle\omega_i,\delta_j\rangle_{1\leq i,j
%\leq p}\right)\cdot m,$$
$$\lambda^p_{{\cal M}}:
{\cal H}{\rm om}_{\HOL}\left(\stackrel{p}{\wedge}\derY,{\cal
M}\right)
\stackrel{\sim}{\to}
{\cal H}{\rm om}_{{\cal V}_0^Y}
\left(\VCERO\otimes_{\HOL}\stackrel{p}{\wedge}\derY,{\cal M}\right),$$
%definido localmente como
%$$\lambda^p_{{\cal
%M}}(\alpha)(P\otimes\delta_1\wedge\cdots\wedge\delta_p)
%= P\cdot\alpha(\delta_1\wedge\cdots\wedge\delta_p).$$
$$\psi^p_{{\cal M}}=\lambda^p_{{\cal M}}\circ\gamma^p_{{\cal M}}:\
\Omega^p_X({\rm log}\ Y)({\cal M})\stackrel{\sim}{\to}
{\cal H}{\rm om}_{\VCERO}
\left(\VCERO\otimes\stackrel{p}{\wedge}\derY,{\cal M}\right).
$$
Locally:
$$\psi^p_{{\cal M}} (\omega_1\wedge\cdots\wedge\omega_p\otimes m)
(P\otimes\delta_1\wedge\cdots\wedge\delta_p) = P\cdot
\det (\langle\omega_i,\delta_j\rangle)_{1\leq i,j
\leq p} \cdot m. $$
%En todos los casos $m$ representa una secci\'{o}n de ${\cal M}$.

\begin{teorema}    \label{TEOREMA}
If ${\cal M}$ is a left $\VCERO$-module  (or, equivalently, is a
$\HOL$-module with an integrable \lo connection),
the complexes of sheaves of $\CC$-vector spaces 
$\Omega^{\bullet}_X(\log Y)({\cal M})$ and $
{\cal H}{\rm om}_{\VCERO}\left({\cal S}p^{\bullet}(\log Y),{\cal M}\right)$
are canonically isomorphic.
\end{teorema}

\begin{prueba}
The general case is solved if we prove the case 
${\cal M}=\VCERO$, using the isomorphisms:
$$\Omega^{\bullet}_X(\log Y)({\cal M}) \cong
\Omega^{\bullet}_X(\log Y)(\VCERO) \otimes_{\VCERO} {\cal M}, $$
$${\cal H}{\rm om}_{\VCERO}\left({\cal S}p^{\bullet}(\log Y),{\cal M} \right)
\cong {\cal H}{\rm om}_{\VCERO}\left({\cal S}p^{\bullet}(\log Y),\VCERO
\right)
\otimes_{\VCERO} {\cal M}.$$
For ${\cal M}=\VCERO$,
 we obtain the right $\VCERO$-isomorphisms 
$$\phi^p = \psi^p_{\VCERO}:\ \Omega^p_X(\log Y)(\VCERO) \to
{\cal H}{\rm om}_{\VCERO}
\left({\cal S}p^{-p}(\log Y),\VCERO \right),$$
whose local expression are:
$$\phi^p \left(
(\omega_1\wedge\cdots\wedge\omega_p)\otimes Q \right)
(P\otimes(\delta_1\wedge
\cdots\wedge\delta_p))=P \cdot \det
\left(\langle\omega_i,\delta_j\rangle\right) \cdot Q .$$
To prove that these isomorphisms produce a 
\underline{isomorphism of complexes} we have to check that they commute
with the differential of the complex.
Thanks to the isomorphism (b) of the proposition \ref{Sp*Y},
$$\D_X[\star Y]\otimes_{\VCERO}{\cal
 S}p^{\bullet}(\log Y)\simeq{\cal S}p^{\bullet}[\star Y],$$
%y considerando que $\DY\otimes_{\VCERO}=\DY$,
we obtain a natural morphism of complexes of sheaves of
right $\VCERO$-modules:
$$\tau^{\bullet}:{\cal H}{\rm om}_{\VCERO}\left({\cal S}p^{\bullet}(\log
Y),\VCERO
  \right) \ \longrightarrow \
  {\cal H}{\rm om}_{{\cal D}_X[\star Y]}\left({\cal S}p^{\bullet}[\star
Y],{\cal
  D}_X[\star Y]\right),$$
locally defined by:
$$\tau^p(\alpha)\left(R\otimes
(\delta_1\wedge\cdots\wedge\delta_p)\right)=
f^{-k}\alpha\left(P\otimes (f\delta_1\wedge\cdots\wedge
f\delta_p)\right)$$ (for any local sections $\alpha$ of
${\cal H}{\rm om}_{\VCERO}\left({\cal S}p^{\bullet}(\log
Y),\VCERO \right)$, $R$ of $\DY$ and
$\delta_1,\cdots,\delta_p$ of $\DER$), where 
$P$ is a local section of $\VCERO$ such that
$ Rf^{-p}=f^{-k}P$ (see lemma \ref{cambiarf}).
The morphisms $\tau^{i}$ are injective, because:
$$\alpha\left(P\otimes (\delta_1\wedge\cdots\wedge
\delta_p)\right)=
\tau^i(\alpha)\left(P\otimes
(\delta_1\wedge\cdots\wedge\delta_p)\right).$$
Let us see the following diagram commutes:
%$$  \Phi_p\circ j_p=\tau_p\circ\phi_p,\ (i\geq 0).$$
$$\begin{array}{ccc}
\Omega^p_X(\log Y)(\VCERO) & \stackrel{j^p}{\longrightarrow} &
\Omega^p_X[\star Y]\left(\DY\right) \\  \ & \ & \ \\
\downarrow {\scriptstyle \phi^p}& \# & \downarrow
{\scriptstyle \Phi^p} \\   \ & \ & \ \\
{\cal H}{\rm om}_{\VCERO}({\cal S}p^p\left(\log Y),\VCERO\right) &
\stackrel{\tau^p}{\longrightarrow} &
{\cal H}{\rm om}_{{\cal D}_X[\star Y]}({\cal S}p^p[\star Y],
{\cal D}_X[\star Y]) \\
\end{array}  $$
for each $p\geq 0$, where the
$\Phi^p$ are the isomorphisms:
$$\Phi^p:\ \Omega^p_X[\star Y]\left({\cal D}_X[\star Y]\right)
\longrightarrow \
   {\cal H}{\rm om}_{{\cal D}_X[\star Y]}\left({\cal
  D}_X[\star
  Y]\otimes\stackrel{p}{\wedge}{\cal D}{\rm er}_{\CC}(\HOL),{\cal
D}_X[\star
  Y]\right),$$
$$\Phi^p((\omega_1\wedge\cdots\wedge\omega_p)\otimes Q)\left(P\otimes
(\delta_1\wedge\cdots
\wedge\delta_p)\right)=P\cdot\det\left(\langle\omega_i
\cdot\delta_j\rangle_{1\leq i,j\leq p}\right)\cdot Q.$$
Given $\omega_1,\cdots,\omega_p$ local sections of $\Omega^1_X(\log
Y)$, $Q$ and $R$ local sections of $\DY$ and
$\delta_1,\cdots,\delta_p$ local sections of $\DER$, we have
$$(\tau^p\circ\phi^p)((\omega_1\wedge\cdots\wedge\omega_p)\otimes Q)
[R\otimes
(\delta_1\cdots
\wedge\delta_p)]=$$ $$f^{-k}\phi_p((\omega_1\wedge\cdots\wedge\omega_p)
\otimes Q)[P
\otimes (f\delta_1\wedge\cdots\wedge f\delta_p)]=$$
$$f^{-k}P\cdot\det
(\langle\omega_i f\delta_j\rangle)\cdot Q =
 R\cdot f^{-p}\det(\langle\omega_i f\delta_j\rangle)\cdot Q =
 R\cdot\det(\langle\omega_i\delta_j\rangle)\cdot Q = $$
$$  \Phi^p\circ j^p
((\omega_1\wedge\cdots\wedge\omega_p)\otimes Q)[R\otimes (\delta_1\wedge
\cdots
\wedge\delta_p)],$$
with $P$ a local section of $\VCERO$ such that $Rf^{-p}=f^{-k}P$.
\vspace{0.2cm}\\
But $\Phi^{\bullet}$, $j^{\bullet}$ and $\tau^{\bullet}$ are 
morphisms of complexes, and $\tau^{\bullet}$ is injective, hence we 
deduce that the $\phi^p$ commute with the differential and so define a
isomorphism of complexes:
$$\phi^{\bullet}:\Omega^{\bullet}_X(\log
Y)\left(\VCERO\right)\longrightarrow
{\cal H}{\rm om}_{\VCERO}
\left({\cal S}p^{\bullet}(\log Y),\VCERO\right),$$
as we wanted to prove.

\end{prueba}

\begin{corolario}
There exists a canonical isomorphism in the derived category:
$$\Omega^{\bullet}_X(\log Y)({\cal M})\cong
{\bf R}{\cal H}{\rm om}_{\VCERO}\left(\HOL,{\cal M}\right).$$
\end{corolario}

\begin{prueba}
By theorem  \ref{splogres},
the complex ${\cal S}p^{\bullet}(\log Y$) is a
locally free resolution of $\HOL$ as left $\VCERO$-module. So, we have 
only to apply the theorem \ref{TEOREMA}.
\end{prueba}

\begin{nota}   \label{rhomdx}
 In the specific case that ${\cal M}=\HOL$,
we have that the complexes
$\Omega^{\bullet}_X(\log Y)$ and $
{\cal H}{\rm om}_{\VCERO}\left({\cal S}p^{\bullet}(\log Y),\HOL\right)$ 
are canonically isomorphic and so, there exists
a canonical isomorphism:
$$\Omega^{\bullet}_X(\log Y)\cong
{\bf R}{\cal H}{\rm om}_{\VCERO}\left(\HOL,\HOL\right).$$
\end{nota}

\begin{nota}
A classical problem is the comparison between the \lo de Rham
complex and the meromorphic de Rham complex relative to a divisor
$Y$,
$$\Omega^{\bullet}_X[\star Y]\cong
{\bf R}{\cal H}{\rm om}_{\DX}\left(\HOL,\OY\right) \cong
{\bf R}{\cal H}{\rm om}_{\VCERO}\left(\HOL,\OY\right) .$$
If $Y$ is a normal crossing divisor, an easy calculation
shows that they are quasi-isomorph (cf. \cite{del_70}).
The same result is true if $Y$ is a strongly
weighted homogeneous free divisor \cite{cas_mon_nar_96}. As a 
consequence of theorem \ref{teorema}, if
$Y$ is an arbitrary free divisor, the meromorphic de Rham complex and
the \lo de Rham complex are quasi-isomorphic if and only if:
$$0 = {\bf R} {\cal H}{\rm om}_{{\cal D}_X}\left({\cal
D}_X\otimes_{\VCERO}^{\bf L} \HOL,
\frac{\HOL[\star Y]}{\HOL}   \right)
\left(= {\bf R} {\cal H}{\rm om}_{\VCERO}\left(\HOL,
\frac{\HOL[\star Y]}{\HOL}  \right) \right) .$$

\end{nota}

\section{Perversity of the \lo complex}

Now we consider the complex
$\DX\otimes_{\VCERO}{\cal S}p^{\bullet}(\log Y)$:
$$0\to\DX\otimes_{{\cal
O}_X}\stackrel{n}{\wedge}\derY\stackrel{{\textstyle \varepsilon}_{-n}}
{\to}\cdots
%$$\cdots\stackrel{{\textstyle \varepsilon}_{-p-1}}{\to}
%\DX\otimes_{{\cal
%O}_X}\stackrel{p}{\wedge}\derY\ \stackrel{{\textstyle
%\varepsilon}_{-p}}
%{\to}\
%\DX\otimes_{\HOL}\stackrel{p-1}
%{\wedge}\derY\stackrel{{\textstyle \varepsilon}_{-p+1}}{\to}\cdots $$
\cdots \stackrel{{\textstyle \varepsilon}_{-2}}
{\to}\DX\otimes_{{\cal
O}_X}\stackrel{1}{\wedge}\derY\stackrel{{\textstyle \varepsilon}_{-1}}
{\to}
\DX,$$
where the local expressions of the morphisms are defined by:
$$\varepsilon_{{\scriptscriptstyle
-p}}(P\otimes(\delta_1\wedge\cdots\wedge\delta_p))
=\sum_{i=1}^p
(-1)^{i-1} P\delta_i\otimes(\delta_1\wedge\cdots\wedge\widehat{\delta_i}
\wedge\cdots\wedge\delta_p) + $$
$$\sum_{1\leq i<j\leq p}
(-1)^{i+j}P\otimes([\delta_i,\delta_j]\wedge\delta_1\wedge\cdots\wedge
\widehat{\delta_i}\wedge\cdots\wedge\widehat{\delta_j}
       \wedge\cdots\wedge\delta_p), \ \ (2\leq p\leq n). $$
$$\varepsilon_{{\scriptscriptstyle -1}} (P\otimes\delta)=P\delta.$$
In the case that $Y$ is a free divisor, we can
work at each point $x$ of $Y$ with a basis $\{\delta_1,\cdots,\delta_n \}$
 of $\derlog$, with $f$ a local reduced equation of $Y$ at $x$.

\begin{proposicion}
If $\{\delta_1,\cdots,\delta_n\}$ is a basis of $\derlog$, and the
sequence $\{\sigma(\delta_1),\cdots,\sigma(\delta_n)\}$ is $\grD$-regular,
it verifies
$$ \sigma\left(\D(\delta_1,\cdots,\delta_n)\right) =
\grD (\sigma(\delta_1),\cdots,\sigma(\delta_n)). $$
\end{proposicion}

\begin{prueba}
The inclusion
$ \grD (\sigma(\delta_1),\cdots,\sigma(\delta_n)) \subset
\sigma\left(\D(\delta_1,\cdots,\delta_n)\right)  $
is clair. Let $F$ be the symbol of an operator $P$
of order $d$,
with $$P=\sum_{i=1}^n P_i\delta_i \in \D(\delta_1,\cdots,\delta_n).$$
We will prove by induction that $F=\sigma(P)$ belongs to
$\grD (\sigma_1,\cdots,\sigma_n),$ with $\sigma_i=\sigma (\delta_i)$.
We will do the induction on the maximum order of the $P_i$
$(i=1,\cdots,n)$, order that we will denote by $k_0$.
 As $P$ has order $d$, $k_0$ is greater or equal to
$d-1$. If $k_0=d-1$, we have:
$$\sigma(P)=\sum_{i\in K}\sigma(P_i)\sigma_i,$$
with $K$ the set of subindexes $j$ such that $P_j$ has
order $k_0$ in $\D$. We suppose that the result holds when $d-1\leq
k_0<m$. Let $F=\sigma (P)$, with $P=\sum_{i=1}^n P_i\delta_i$ and
$k_0 = m$. There are two possibilities:
\begin{enumerate}
\item $F=\sigma(P)=\sum_{i\in K} \sigma(P_i)\sigma_i\in
\grD(\sigma_1,\cdots,\sigma_n)$, as we wanted to prove.
\item $ \sum_{i\in K} \sigma(P_i)\sigma_i = 0 $.
\end{enumerate}
In this last case, as
$\{\sigma_1,\cdots,\sigma_n\}$ is a 
$\grD$-regular sequence, if we call $F_i$ the symbol $\sigma(P_i)$
 in the case that
$i\in K$ and 0 otherwise,  we have:
$$\left( F_1,\cdots,F_n\right) =
\sum_{i<j} F_{ij} (0,\cdots,0,\stackrel{\stackrel{i}{\smile}}{\sigma_j},
0,\cdots,0,\stackrel{\stackrel{j}{\smile}}{- \sigma_i},0,\cdots,0), $$
with $F_{ij}\in \grD$  homogeneous polynomials of order $m-1$.
We choose, for $1\leq i<j \leq n$,
operators $Q_{ij}$, of order $m-1$ in
$\D$, such that
$\sigma (Q_{ij})=F_{ij}$, and define:
$$\left( Q_1,\cdots,Q_n\right) = \left(
P_1,\cdots,P_n\right) -
\sum_{i<j} Q_{ij} \left((0,\cdots,0,\stackrel{\stackrel{i}{\smile}}{\delta_j},
0,\cdots,0,\stackrel{\stackrel{j}{\smile}}{- \delta_i},0,\cdots,0) -
{\underline{\alpha}}_{ij}
\right),$$
where $\underline{\alpha}_{ij}$ are the
vectors with $n$ coordinates in $\hol$ defined by the relations:
$$[\delta_i,\delta_j] = \sum_{k=1}^n a_{ij}^k \delta_k =
\underline{\alpha}_{ij}\bullet\underline{\delta},$$
with $\underline{\delta}=(\delta_1,\cdots,\delta_n)$.
These $Q_i$, of order $m$ in $\D$, verify
$$\left( \sigma_m(Q_1),\cdots,\sigma_m(Q_n)\right) = $$
$$\left( F_1,\cdots,F_n\right) -
\sum_{i<j} F_{ij} (0,\cdots,0,\stackrel{\stackrel{i}{\smile}}{\sigma_j},
0,\cdots,0,\stackrel{\stackrel{j}{\smile}}{- \sigma_i},0,\cdots,0)=0. $$
So, $Q_i$ has order $m-1$ in $\D$. Moreover,
$$ \sum_{i=1}^n Q_i \delta_i =
\sum_{i=1}^n P_i\delta_i -
\sum_{i<j} Q_{ij} \left( \delta_i\delta_j - \delta_j\delta_i -
[\delta_i,\delta_j] \right) = \sum_{i=1}^n P_i\delta_i = P .$$
We apply the induction hypothesis to $F=\sigma(P)$, with $ P=
\sum_{i=1}^n Q_i \delta_i $, and obtain:
$$\sigma(P)\in \grD(\sigma_1,\cdots,\sigma_n).$$
\end{prueba}

\begin{proposicion}
Let $\{\delta_1,\cdots,\delta_n\}$ be a basis of $\derlog$.
If the sequence $\sigma(\delta_1),\cdots,\sigma(\delta_n)$ is a
$\grD$-regular sequence in $\grD$, the complex
$\D\otimes_{{\cal V}_0^f({\cal D})}{\cal S}p^{\bullet}({\rm log}\ f)$ is a
resolution of the quotient module $\frac{\D}{\D(\delta_1,\cdots,\delta_n)}$.
\end{proposicion}

\begin{prueba}
We consider the complex $\D\otimes_{{\cal V}_0^f({\cal D})}{\cal S}
p^{\bullet}({\rm log}
\ f)$. We can augment this complex of $\D$-modules by another morphism:
$$\varepsilon_0:\D\to \frac{\D}{\D(\delta_1,\cdots,\delta_n)} ,\ \
\varepsilon_0(P)=P + \D(\delta_1,\cdots,\delta_n).$$
We denote by $\D\otimes_{{\cal V}_0^f({\cal D})}\widetilde{{\cal S}}
p^{\bullet}(\log f)$ the new complex. To prove that this new complex 
is exact, we define a discrete filtration 
$G^{\bullet}$ 
such that the graded complex be exact
 (cf. \cite[lemma 3.16]{bjo_79}):
$$G^k\left(\D\ \otimes_{\hol}\stackrel{p}{\wedge}\derlog\right)
=F^{k-p}\left(\D\right)\otimes_{\hol}\stackrel{p}{\wedge}\derlog,$$
$$G^k\left(\frac{\D}{\D(\delta_1,\cdots,\delta_n)}\right) =
\frac{F^k(\D) + \D\cdot(\delta_1,\cdots,\delta_n) }
{\D (\delta_1,\cdots,\delta_n)}.$$
Clairly the filtration is compatible with the differential of the complex.
%$$\varepsilon_{-p}\left(
%G^k\left(\D\ \otimes_{\hol}\stackrel{p}{\wedge}\derlog\right)  \right)
%\subset G^k\left(\D\
%\otimes_{\hol}\stackrel{p-1}{\wedge}\derlog\right),
%$$
%$$ \varepsilon_0 \left( G^k(\D)\right)\subset
%G^k\left(\frac{\D}{\D(\delta_1,\cdots,\delta_n)}\right) $$
%donde $F^{\bullet}\left(\vcero\right)$ es la filtraci\'{o}n por el
%grado de los operadores.\\
Moreover:
$${\rm
G}{\rm r}_{G^{\bu}}\left(\D\ \otimes\stackrel{p}{\wedge}\derlog\right)=
\grD[-p]\otimes\stackrel{p}{\wedge}\derlog,$$
and, by the previous proposition,
$${\rm G}{\rm r}_{G^{\bu}}\left(\frac{\D}{\D(\delta_1,\cdots,\delta_n)}\right)
= \frac{\grD}{\sigma\left(\D\cdot(\delta_1,\cdots,\delta_n)\right)} =
 \frac{\grD}{\grD\cdot(\sigma(\delta_1),\cdots,\sigma(\delta_n) )}.$$
We consider the complex
${\rm G}{\rm r}_{G^{\bu}}\left(\D\otimes_{{\cal V}_0^f({\cal D})}
\widetilde{{\cal
S}}p^{\bullet}(\log\ f)\right):$
$$0\to\grD[-n]\otimes_{{\cal
O}}\stackrel{n}{\wedge}\derlog\stackrel{\psi_{-n}}{\to}\cdots
%$$\cdots\stackrel{\psi_{-p-1}}{\to}
%\grD[-p]\otimes_{{\cal O}}\stackrel{p}{\wedge}\derlog\
%\stackrel{\psi_{-p}}{\to}$$
%$$\stackrel{\psi_{-p}}{\to}\ \grD[-p+1]\otimes_{{\cal
%O}}\stackrel{p-1}{\wedge}\derlog\stackrel{\psi_{-p+1}}{\to}\cdots $$
\stackrel{\psi_{-2}}{\to}\grD[-1]\otimes_{{\cal
O}}\stackrel{1}{\wedge}\derlog $$ $$\stackrel{\psi_{-1}}{\to}
\grD\stackrel{\psi_0}{\to}
\frac{\grD}{\grD\cdot(\sigma(\delta_1),\cdots,\sigma(\delta_n))}
\to 0,$$
where the local expression of the differential is defined by:
$$\psi_{-p}(G\otimes\delta_{j_1}\wedge\cdots\wedge\delta_{j_p})=
\sum_{i=1}^p
(-1)^{i-1}
G\sigma(\delta_{j_i})\otimes\delta_{j_1}\wedge\cdots\wedge
\widehat{\delta_{j_i}}
\wedge\cdots\wedge\delta_{j_p},\ \  (2\leq p\leq n), $$
$$\psi_{-1}(G\otimes\delta_i)=G\sigma(\delta_i),$$ $$\psi_0(G)=G +
\grD\cdot(\sigma(\delta_1),\cdots,\sigma(\delta_n)). $$
%Siendo $ \{ \delta_1,\cdots,\delta_n\}$ una base (local) de $\derlog$.
This complex is the Koszul complex of the ring
$\grD$ with respect to the sequence
 $\sigma(\delta_1),\cdots,\sigma(\delta_n)$.
So we deduce that, if
the sequence $\sigma(\delta_1),\cdots,\sigma(\delta_n)$ is
$\grD$-regular in $\grD$,  the complex
$${\rm G}{\rm r}_{G^{\bu}}\left(\D\otimes_{{\cal V}_0^f({\cal D})}
\widetilde{{\cal
S}}p^{\bullet}(\logf)\right)$$ is exact. So,
the complex
$\D\otimes_{{\cal V}_0^f({\cal D})}\widetilde{{\cal
S}}p^{\bullet}(\logf)$ is exact too, and
$\D\otimes_{{\cal V}_0^f({\cal D})}{\cal S}p^{\bullet}{\rm(log} \ f)$ is a 
resolution
of  $\frac{\D}{\D(\delta_1,\cdots,\delta_n)}$.
\end{prueba}

\begin{corolario}         \label{perversion}
Let $Y$ be a free divisor. With the conditions of the previous proposition
(for each point $x$ of $Y$, there exists a basis
$\{\delta_1,\cdots,\delta_n\}$ of $\derlog$ such that
the sequence $\sigma(\delta_1),\cdots,\sigma(\delta_n)$ is a
$\grD$-regular sequence), the sheaf
$\Omega^{\bullet}_X(\logY)$ is a perverse sheaf.
\end{corolario}

\begin{prueba}
With the same conditions of the previous proposition,
 the homology of the complex
$\DX\otimes_{\VCERO}{\cal S}p^{\bullet}{\rm (log} \ Y)$ is concentrated 
in degree 0. All its homology groups are zero except the group
in degree 0, which
verifies:
$$h^0\left(\DX\otimes_{\VCERO}{\cal S}p^{\bullet}{\rm (log}
\ Y)\right) =
\frac{\DX}{\DX\cdot \derY} =
\frac{{\cal D}_X}{{\cal D}_X\cdot(\delta_1,\cdots,\delta_n)}={\cal E},$$
where $\{\delta_1,\cdots,\delta_n\}$ is a local basis of $\derY$.
 But ${\cal E}$
is a holonomic $\DX$-module because:
$${\cal G}{\rm r}_F({\cal E}) = \frac{\GrD}
{\left(\sigma(\delta_1),\cdots,\sigma(\delta_n)\right)}$$
has dimension $n$ (using the fact that
 $\sigma(\delta_1),\cdots,\sigma(\delta_n)$
is a $\GrD$-regular sequence).
 So (using remark
\ref{rhomdx} for the first equality and 
teorema \ref{splogres} for the last equality)):
$$\Omega_X^{\bullet}({\rm log}\ Y) = {\bf R} {\cal H}{\rm om}_{{\cal
D}_X}\left({\cal D}_X\otimes_{\VCERO}^{\bf L} \HOL,
\HOL   \right) = $$
$${\bf R} {\cal H}{\rm om}_{{\cal
D}_X}\left({\cal D}_X\otimes_{\VCERO}{\cal S}p^{\bullet}{\rm (log}
\ Y), \HOL   \right) =
 {\bf R} {\cal H}{\rm om}_{{\cal D}_X}\left(
\frac{{\cal D}_X}{{\cal D}_X(\delta_1,\cdots,\delta_n)} ,
 \HOL   \right) $$
is a perverse sheaf (as solution of a holonomic $\DX$-module, 
cf. \cite{meb_formalisme}).

%y puesto que el complejo ${\cal S}p^{\bullet}$(log $Y$) es una
%resoluci\'{o}n localmente libre de
%$\HOL$ como $\VCERO$-m\'{o}dulo a la izquierda (teorema
%\ref{splogres}),
%obtenemos el resultado que buscamos.
\end{prueba}

\begin{corolario}
Let $Y$ be any divisor in $X$, with $\dim_{\CC} X=2$.
Then $\Omega^{\bullet}_X($log $Y)$ is a perverse sheaf.
\end{corolario}

\begin{prueba}
 We know that, if $\dim_{\CC} X=2$, any divisor $Y$ in $X$ is free
\cite{ksaito_log}. So, we have only to check that the other hypothesis
of the previous corollary holds. We consider the symbols
$\{\sigma_1,\sigma_2\}$ of a basis $\{\delta_1,\delta_2\}$
of $\derlog$, where
$f$ is a reduced equation
of $Y$. We have to see that they form a 
$\grD$-regular sequence.
If they do not, they have a common factor
$g\in \hol$, because they are symbols of operators of order 1. If $g$ is
a unit, we divide one of them by $g$ and eliminate the common factor.
 If $g$ is not a unit, it would be in contradiction with Saito's Criterion,
because the determinant of 
the coefficients of the basis
$\{\delta_1,\delta_2 \}$ would have as factor $g^2$, with $g$ not invertible, 
and this determinant has to be equal to $f$ multiplied by a unit.
\end{prueba}

\begin{nota}
The regularity of the sequence of the symbols of a basis of $\derlog$ in
$\grD$ is not necessary for the perversity of the \lo
de Rham complex. For example, if $X=\CC^3$ and $Y\equiv\{f=0\}$,
with
$f=xy(x+y)(y+tx)$, $f$ is a free divisor such that the graded complex
$$ {\cal G}{\rm r}_{G^{\bu}}(
{\cal D}_X\otimes_{\VCERO} {\cal S}p^{\bullet}(\logY))=
K(\sigma(\delta_1),\sigma(\delta_1),\sigma(\delta_3);{\cal G}
{\rm r}_{F^{\bu}}({\cal D}_X))$$
is not concentrated in degree 0, but the complex
$${\cal D}_X\otimes_{{\cal V}_0^Y({\cal D}_X)}
{\cal S}p^{\bullet}(\logY)$$
is. Moreover, in this case the dimension of
$\frac{{\cal D}_X}{{\cal D}_X\cdot(\delta_1,\delta_2,\delta_3)}$ is 3
and so, $\Omega^{\bullet}_X(\log Y)$ is a perverse sheaf.

\end{nota}
\bigskip

\medskip
{\small
\bibliographystyle{plain}

}

\end{document}